\documentclass[12pt]{amsart}
\usepackage{amssymb}
\usepackage[dvips]{graphics}
\ifx\mathrm\undefined\let\mathrm\rm\fi
\ifx\mathbf\undefined\let\mathbf\bf\fi
\ifx\mathfrak\undefined\let\mathfrak\frak\fi
\ifx\mathcal\undefined\let\mathcal\cal\fi
\ifx\mathbb\undefined\let\mathbb\Bbb\fi
\ifx\emph\undefined\let\emph\it\fi
 at9.98pt

\newcommand{\Z}{{\mathbb Z}}

\newcommand{\R}{{\mathbb R}}
\newcommand{\C}{{\mathbb C}}
\newcommand{\Q}{{\mathbb Q}}

\newcommand{\dontprint}[1]
{\relax}

\newtheorem%
{thm}{Theorem}[section]
\newtheorem%
{proposition}[thm]{Proposition}
\newtheorem%
{lemma}[thm]{Lemma}
\newtheorem%
{lemmadef}[thm]{Lemma-Definition}
\newtheorem%
{corollary}[thm]{Corollary}
\newtheorem%
{conjecture}[thm]{Conjecture}

\newcommand{\bea}{\begin{eqnarray*}}
\newcommand{\eea}{\end{eqnarray*}}
\newcommand{\bean}{\begin{eqnarray}}
\newcommand{\eean}{\end{eqnarray}}

\newcommand{\nc}{\newcommand}
\nc{\on}{\operatorname}
\nc{\al}{\alpha}
\nc{\ri}{\rangle}
\nc{\lef}{\langle}
\nc{\W}{{\mathcal W}}
\nc{\La}{\Lambda}
\nc{\ep}{\epsilon}
\nc{\Om}{\Omega}
\newcommand{\be}{\begin{displaymath}}
\newcommand{\ee}{\end{displaymath}}

\nc{\PCr}{{ \Bbb P  (\C[x])^r   }}

\newcommand{\Zo}{\Z_{\rm{odd}}\,}

\begin{document}

\title[Finite Order Invariants for $(n,2)$-Torus Knots]
{Finite Order Invariants for $(n,2)$-Torus Knots \\ and the Curve
$Y^2 = X^3 + X^2$}
\author[S. Tyurina and A. Varchenko]
{Svetlana Tyurina${}^{*}$ \ {}\ and \ {} \ Alexander Varchenko${}^{**,1}$ }
\thanks{${}^1$Supported in part by NSF grant  DMS-0244579}

\begin{abstract}
We describe the algebra of finite order invariants on the set of all
$(n,2)$-torus knots.

\end{abstract}

\maketitle

\centerline{\it ${}^{*}$OASIS,
University of North Carolina at Chapel Hill,}
\centerline{\it Chapel Hill, NC 27599, USA}
\centerline{svetlana@email.unc.edu}

\medskip
\centerline{\it ${}^{**}$Department of Mathematics,
University of North Carolina at Chapel Hill,}
\centerline{\it Chapel Hill, NC 27599-3250, USA}
\centerline{anv@email.unc.edu}

\bigskip

\centerline{February, 2004}

\bigskip

This paper is an extended exposition of the talk \cite{Ty} given by the
first author. The authors thank S. Duzhin and A. Sossinsky for interest to this
work and S. Chmutov for useful  discussions.

\bigskip


Consider the $\Q$-algebra $V$ of Vassiliev finite order knot invariants, see
for example \cite{B, CDL}. The algebra is filtered,
\bea
V_0 \subset V_1 \subset \dots \subset V_k \subset \dots \subset V ,
\eea
the vector subspace $V_k \subset V$ consists of  knot invariants of order
not greater than $k$. We have $V_k \cdot V_l \subset V_{k+l}$. 

The subspace $V_0$ is of
dimension $1$ and consists of invariants taking the same value
on all knots. It is known that $V_1 = V_0$, \ $\dim V_2/V_1 = 
\dim V_3/V_2 = 1$. The generator in $V_2/V_1$ is given by the knot
invariant $x$ of order $2$ which takes value $0$ on the trivial knot and
value $8$ on the trefold. 
The generator in $V_3/V_2$ is given by the knot
invariant $y$ of order $3$ which takes value $0$ on the trivial knot, takes
value $24$ on the trefold, and takes value $-24$ on its  mirror image. 
Those conditions
determine $x$ and $y$ uniquely,  see for example \cite{L}.

It is known that the  space  $V_k$ has finite dimension 
fast growing with  $k$, see  for example \cite{CD, D, Z}.

By definition the algebra $V$ is an
algebra of  certain special functions on the set $K$ of
all knots in $\R^3$ considered up to isotopy.
 
Let $T \subset K$ be the subset of toric knots of type $(n,2)$,
$n = \pm 1, \pm 3, \dots$.
Here $(1,2)$ and $(-1,2)$  
denote the trivial knot, 
$(3,2)$ is the trefoil, $(-3,2)$  its mirror image,  and so on.

Consider the algebra $A$ of functions on $T$, which is the restriction of
$V$ to $T$, i.e.  $A = V|_T$. The algebra $A$ is filtered,
\bea
A_0 \subset A_1 \subset \dots \subset A_k \subset \dots \subset A ,
\eea
where $A_k = V_k|_{T}$ for any $k$. Our goal is to describe $A$.

Let $X \in A_2$ and $Y \in A_3$ be the image of $x$ and $y$, 
respectively, under the natural projection $V \to A$.

\bigskip

\noindent
{\bf Theorem.}
{\it The algebra $A$ is generated by $X$ and $Y$ and is isomorphic to
the algebra $\Q[X,Y]/(X^3 + X^2 - Y^2)\Q[X,Y]$,\ where
$(X^3 + X^2 - Y^2)\Q[X,Y]\ \subset \ \Q[X,Y]$ 
is the ideal generated by the polynomial $X^3 + X^2 - Y^2$. 
We have $\dim A_0 = 1$, $\dim A_1/A_0 = 0$,
$\dim A_k/A_{k-1} = 1$ for $k > 1$. 
The generator in $ A_{2l}/A_{2l-1}$ is given by $X^l$
and the 
generator in $ A_{2l+1}/A_{2l}$ is given by $X^{l-1}Y$
for all $l > 0$. }

\bigskip

{\it Proof.} Denote by $\Zo$ the set of all odd integers.
For   $n\in \Zo$ denote by $[n]$ the  
torus knot of type $(n,2)$. 

An element $f \in A$ defines  a function
$ \Zo \to \Q$, $ n \mapsto f( [n] )$, and is uniquely
determined by that function. Thus $A$ can be considered
as an algebra of certain functions on $\Zo$.

\bigskip

\noindent
{\bf Lemma. }
{\it
\begin{enumerate}
\item[$\bullet$]
 If $f : \Zo \to \Q$ belongs to $A_k$ for some $k$, 
then $f$ is a polynomial of degree not greater than $k$.
\item[$\bullet$]
 If $f : \Zo \to \Q$ belongs to $A$, then $f(1) = f(-1)$.
\end{enumerate}
}
\hfill
$\square$

\bigskip

The lemma is a direct corollary of definitions.

\bigskip

We have $X : \Zo \to \Q,\ n \mapsto n^2-1$, and 
$Y : \Zo \to \Q, \ n \mapsto n^3-n$. This  gives the relation
$Y^2 = X^3 + X^2$. 

It is easy to see that 
all polynomials $ f : \Zo \to \Q$ with  property
$f(1) = f(-1)$ are linear combinations of monomials $X^l$ and $X^{l-1}Y$
of degree $2l$ and $2l+1$ respectively. 

The theorem is proved.
\hfill
$\square$

\bigskip

\noindent
{\bf Remarks.}
\begin{enumerate}
\item[$\bullet$]
After this paper had been written, S. Chmutov informed us about 
paper \cite{T}, where R. Trapp in particular shows that any element 
$f \in A$ is a polynomial function on $\Zo$ and $f$ can be written as a polynomial
in $X$ and $Y$.

\item[$\bullet$]  
S. Chmutov informed us that the shapes, similar to the shape of
our curve $Y^2 = X^3 + X^2$, appeared in \cite{W}, where
S. Willerton discusses statistics of points $(x(k),y(k)) \in \Q^2$
for arbitrary knots $k$.

\item[$\bullet$] According to our theorem the algebra $V|_T$ is isomorphic to the
algebra of regular functions on the affine curve $Y^2 = X^3 + X^2$.
One may wander what kind of algebraic varieties one obtains considering restrictions
of $V$ to other reasonable subsets of the set of all knots.

\end{enumerate}

\end{document}